\documentclass{birkmult}
 \theoremstyle{definition}
 
 \theoremstyle{remark}

 \numberwithin{equation}{section}
\newcommand{\fer}[1]{(\ref{#1})}

\newcommand{\R}{\mathbb{R}}

\newcommand {\al} {\alpha}

\newcommand {\lb} {\lambda}

\newcommand {\Chi} {{\bf \raise 2pt \hbox{$\chi$}} }

\newcommand {\cac} { {\mathcal C} }

\newcommand {\f}   {\frac}
\newcommand {\p}   {\partial}
\newcommand{\beq}{\begin{equation}}
\newcommand{\beqa}{\begin{eqnarray}}
\newcommand{\bea} {\begin{array}{ll}}
\newcommand{\beqan}{\begin{eqnarray*}}
\newcommand{\eeq}{\end{equation}}
\newcommand{\eeqa}{\end{eqnarray}}
\newcommand{\eeqan}{\end{eqnarray*}}
\newcommand{\eea} {\end{array}}

\newcommand{\ep}{\epsilon}

\usepackage{setspace}
\usepackage{hyperref}
\begin{document}

%
%
%
%
%
%
%
%
%
\title[Calibration of a Population Model]
 {On the Calibration of a Size-Structured Population Model from
Experimental Data}
\author[Doumic]{Marie Doumic}
\address{
INRIA Rocquencourt, projet BANG,
Domaine de Voluceau, BP 105, 78153 Rocquencourt, France
}

\email{marie.doumic@inria.fr}

\author[Maia]{Pedro Maia}
\address{UFRJ, Cidade Universitaria, Caixa-Postal: 68530, Rio de Janeiro, RJ 21945-970, Brazil.}

\email{pd\_maia@hotmail.com}


\author[Zubelli]{Jorge P. Zubelli}

\address{%
IMPA, 
Est. D. Castorina 110, 
Rio de Janeiro, RJ 22460-320, 
Brazil.
}

\email{zubelli@impa.br}
\thanks{Work supported by CNPq Grant \#302161/2003-1 and 474085/2003-1 , CNPq-INRIA agreement, Brazil-France project.}



\date{January 1, 2009}

\begin{abstract}
The aim of this work is twofold. First, we survey the techniques developed in \cite{PZ,DPZ} to reconstruct the division (birth) rate from the cell volume distribution
data in certain structured population models. Secondly, we implement such techniques on  experimental cell volume distributions available in the
literature so as to validate the theoretical and numerical results. As a proof of
concept, we use the data reported in the {\em classical} work of Kubitschek~\cite{K69} concerning Escherichia coli {\em in vitro} experiments measured by means of a  Coulter transducer-multichannel analyzer system (Coulter Electronics, Inc., Hialeah, Fla, USA.)
Despite the rather old measurement technology, the reconstructed division rates still 
display potentially useful biological features.
%

\end{abstract}
%

\maketitle
\section{Introduction}

Structured growth models have been the subject of much attention for the last 
decades \cite{MD,BP}. Amongst such models, equal division ones play a crucial
role for their relative simplicity. In equal division models, each cell of 
volume $2x$ subdivides into two cells of size $x$ according to a certain rate $B(2x)$.
In order to understand better such models and related microbiological mechanisms, it would be important to calibrate the birth rate $B(\cdot)$ by means of observed cell densities
$n(t,x)$, where $t$ is the time parameter. 
Thus, in this article, we consider the problem of calibrating from measured data the following size-structured model for cell division:
\begin{equation}\label{eq1}
\left \{ \begin{array}{l}
\f{\p}{\p t} n(t,x) +	\f{\p}{\p x}[g(x) n(t,x)]  + B(x) n(t,x)  = 4  B(2 x) n(t,2 x), \quad x, t \geq 0,   \\
g(x=0) n(t,x=0)=0, t> 0,   \\
n(0,x) = n^{0}(x) \ge 0 .
\end{array} \right.
\end{equation} 
The function $g$ describes the microscopic growth behaviour of individual cells.

Usually two different models have been considered. The first one consists of constant $g$, which without loss of generality we take to be one. This is usually called {\em linear growth model} and leads to characteristic lines of the form $x(t)=x_0 + t$.
The second one is $g(x) = \kappa x$, and leads to characteristic lines
of the form $x(t) = x_0 \exp(\kappa t)$ where $\kappa$ is a constant. It is called
{\em exponential growth model}. 
Obviously both of them subsume some kind of unlimited growth of individual organisms.
There has been substantial amount of debate in the literature concerning the 
correct way of modeling the function $g$. See for example \cite{Cooper,Koch,Mitchison} and references therein.

Two macroscopic quantities of biological interest are naturally computed from the model~(\ref{eq1}). 
The total cell quantity  
$\mathcal{N}(t)= \int_0^\infty n(t,x) dx $
and  the total biomass 
$\mathcal{M}(t)=\int_0^\infty x n(t,x) dx  \mbox{ .}$
Integrating equation~(\ref{eq1}) yields
 $\dot{\mathcal N}(t)=\int_0^\infty B(x) n(t,x)dx \mbox{ .}$
This  means that number of cells increase only by division.
Integrating the identity function times~(\ref{eq1}) yields
$\dot{\mathcal M}(t)=\int_0^\infty g(x) n(t,x)dx \mbox{ .}$
This  means that the biomass increases only by nutrient uptake.

%
%

Under reasonable assumptions
on $B$, the long time behavior of $n(t,x)$ is of the form
$n(t,x) \approx m_0 N(x) e^{\lambda_0 t} $. More precisely, it is proved in \cite{GDP2009} (see also \cite{PR, M} for the case $g(x)=1$) that, under fairly general conditions on the coefficients, 
there is a unique solution $(N,\lb_0)$  to the following eigenvalue problem
\beq
\label{eq:celldiv1}
\left \{ \begin{array}{l}
 \f{\p}{\p x} (g N) + (\lb_0+B(x)) N =4  B(2 x) N(2 x), \qquad x > 0,
\\
g(x=0)N(x=0)=0,
\\
N(x)>0 \; \text{ for } x>0, \qquad \int_0^\infty N(x)dx =1, 
\end{array} \right.
\eeq
where $\lb_0>0$ and $x^\alpha g N \in L^\infty \cap L^1$ for all $\alpha \geq 0$.
Furthermore, it was also shown in \cite{GDP2009, PR,MMP} that
\begin{equation} \label{eq14p}
n(t,x) e^{-\lb_0 t} {\;}_{\overrightarrow{\; t \rightarrow \infty \;}}\;  m_0 N(x), \quad \mbox{ in }\;L^1(\R_+,\phi(x)dx),
\end{equation}
where the weight $\phi$ is the unique solution to the adjoint problem to Equation~(\ref{eq:celldiv1}).
%
In other words, $\lb_0$ is the growth rate of such a system and is usually called ``Malthus parameter'' in population biology. 
By integrating the equation and integrating it against the weight $x,$ we also know  that $\lb_0$ is related to $N$ by the relation
\beq\label{eq:lbn}
\lb_{0}= {\int_0^{\infty} g(x) N (x) dx}\bigg{/}{\int_0^{\infty} x N dx}.
\eeq


%

The plan for the article is the following: In Section~\ref{bio}, we discuss
some biological preliminaries of the model under consideration. 
In Section~\ref{method} we describe the
theoretical aspects of the inverse problem under consideration. This will give us the
underlying theoretical basis for the calibration techniques. In Section~\ref{main}, we present the calibration results. In Section~\ref{disc} we conclude with a discussion
of the results and suggestions for further research.

\section{Biological Preliminaries} \label{bio}

In many situations it is important to understand when an organism
is mature and ready to reproduce. The concept of size at
maturation was developed as a theoretical tool to express how
individuals deal with environmental variability within growth
rates. 
%
%
%
In addition, individuals may have unequal fitness and might
respond differently to the same environmental conditions. In this
context, we can interpret the bacterial-growth experiments of
Kubitschek \cite{K69} as a very special situation: (1) the cells grow in a
chemostat and all individuals have the same amount of nutrients,
in all generations of the population. (2) the cells are clones
from each other and have the same fitness, responding equally to
the same environmental effects. (3) Kubitschek \cite{K69} also did a
synchronic experiment, using temperature control techniques, so
that all the cells would divide at the same time.

One might argue that in this situation, with no difference in the
fitness of the individuals and in the environmental conditions, 
there should be almost no variability in the preferred size of
division. However, the reconstructions of $B(x)$ detect an
intrinsic variability that is inherent from the division process.
This justifies the importance of the reconstructed division rate.
Moreover, \textit{E. coli} is commonly used as a model for more
complex organisms and the present work can be considered a first
step toward the modelling of more complex cell division phenomena.
%

We now address the biological suitability of the mathematical model.
The specific size-structured model under consideration in this work is suitable to microorganisms with the following characteristics:
\begin{quote}
(a) All the cells of a  given class of microorganisms grow
according to some deterministic rule.\\
(b) All cells, after attaining a sufficient size, divide into two,
and only two, identical daughter cells.\\
(c) The cell's age (generation) does not influence the growth law
and the preferable size of division.\\
(d) The total number of individuals of the population increases
exponentially.\\
(e) There is a negligible number of deaths during normal growth.\\
\end{quote}

We will now argue that the statements above are reasonable for the
bacteria \textit{E. Coli} growing under certain conditions: 

Assumption (a) is in the spirit of \cite{Cooper}, where the author argues for the existence of a growth law of cells that can be found
and understood. 

Perhaps the most delicate assumption is (b). The cell division
process may result in daughters with unequal sizes. Equal division
does not occur surely for all individuals of a microbiological
population, but we may quantify how much they deviate from a
perfect division by a statistical description, as explained in the
work of Koch~\cite{Koch}. 

The bacteria \textit{E. Coli} is known for
its quite small variation in its subdivisions \cite{Koch,Trueba,Marr}.
%
It is usually assumed (see \cite{Koch}) that the fluctuations in
the critical size-at-division of individual cells is random and
not correlated with other cell cycle events in the current cell
generation or in earlier generations. This justifies (c) for the
bacteria \textit{E. Coli}.

It is well known that \textit{E. coli} cells in normal growth conditions
achieves an exponential phase \cite{Prescott} thus confirming claim (d)
in agreement with Equation~\ref{eq14p}.

Finally, for statement~(e), cell death during normal growth is
considered by Koch \cite{Koch} as a minor source of
variability. It is true that some cells may be dead or moribund
and may distort the distribution if they are clustered at some
cell size range. But these effects can be neglected if there is
plenty of nutrient, which occurs in the exponential phase of the
bacterial growth. This hypothesis must be carefully revised when
the population is submitted to starvation and other stressful
types of experiments.
The present model would not apply for instance, to the budding
yeast \textit{S. cerevisae} \cite{Porro}. 
However, it could easily be adapted to take into account a known death rate.

\section{Methodology: Inverse problems and theoretical results} \label{method} 

Based on the limiting behavior of the distribution $n(t,x)$ as $t$ grows, the problem
of calibration of $B$ assuming $g\equiv 1$ reduces to the following: How can we estimate the division rate $B$ from the knowledge of the steady dynamics $N$ and  $\lb_0$ ? 
This corresponds to finding $B$ a solution to 
\beq\label{eq:exact:inverse}
 4B(2x)N(2x) - B(x) N(x) =L(x):= \f{\p}{\p x} N(x) + \lb_0 N(x) ,\quad x >  0,
\eeq
assuming that $(N,\lb_0)$ is known, or, thanks to \fer{eq:lbn}, that $N$ is known. 
%
However, in practical applications we have only an \emph{approximate} knowledge of $(N,\lb_0),$ given by noisy data $(N_\ep,\lb_\ep),$ with $N_\ep\in L^2_+(\R_+)$ for instance.
This means that we have no way of controlling $\f{\p}{\p x} N_\ep,$ so we cannot control the precision of a solution $B_\ep$ to problem (\ref{eq:exact:inverse}) when a perturbed $N_\ep$ replaces  $N$. Thus, in the context of noisy data, the inverse problem under consideration is ill-posed and  regularization is needed.

As alluded above, $ N_\ep $ denotes the measured stable distribution. The subindex refers to 
the difference in the appropriate norm to the value $N$ associated to the true
$B$ of the model. The precise value of $\ep$ is obviously unknown, but it depends 
on a number of factors such as accuracy of the measurements, quality of the model
and proximity of the interpolation leading from a finite amount of data to the function
with domain $\mathbb{R}_{+}$.  

Following~\cite{DPZ} three methods will be used, namely: quasi-reversibility,  
 filtering, and a hybrid of these.

\paragraph{Regularization by Quasi-Reversibility}

To regularize the problem by the quasi-reversibility proposed in \cite{PZ} we work with 
\beq
\left\{
\begin{array}{lll}
\al \f{\p}{\p y} (B_{\ep,\al} N_\ep) & + & 4B_{\ep,\al}(y) N_\ep(y)   =  B_{\ep,\al}\big(\f y 2\big) N_\ep\big(\f y 2\big)+ \lb_{\ep,\alpha} N_\ep\big(\f y 2\big)+2 \f{\p}{\p y}\biggl( N_\ep\big(\f y 2\big)\biggr), \qquad y>0,
\\
(B_{\ep,\al} N_\ep)(0) & = & 0.
\end{array} \right.
\label{eq:invfull}
\eeq
According to the eigenvalue theory \cite{DPZ}, we have to choose 
\beq\label{eq:conserv2:PZ}
\lb_{\ep,\alpha}= \left({\int_0^{\infty} N_\ep dx}\right)\bigg{/}\left({\int_0^{\infty} x N_\ep dx + \f{\alpha}{4}\int_0^{\infty} N_\ep dx}\right).
\eeq
The method thus consists in solving numericaly the Equation~\fer{eq:invfull}. This was analyzed in detail in Section~3 of \cite{DPZ}. 

\paragraph{Regularization by Filtering} 
In this method we basicaly filter the data so as to smooth out the noisy component of
the measurements. 
For $\alpha>0$,  we introduce
$$
\rho_\alpha(x) = \f{1}{\alpha}\rho(\f{x}{\alpha}),\qquad \rho\in \cac_c^\infty(\R),\quad\int_0^{\infty}\rho(x)\,dx=1,\quad\rho\geq  0,\quad \mathrm{Supp}(\rho)\subset [0,1], 
$$
and we replace in (\ref{eq:exact:inverse}) the term $\f{\p}{\p x} N_\ep + \lb_0 N_\ep$ by
$$\bigl(\f{\p}{\p x} N_\ep + \lb_{\ep,\alpha} N_\ep\bigr)* \rho_\alpha (x)=
N_\ep * \bigl(\f{\p}{\p x} \rho_\alpha + \lb_{\ep,\alpha} \rho_\alpha \bigr)(x)=\int_0^{\infty} N_\ep (x') \bigl(\f{\p}{\p x} \rho_\alpha +\lb_{\ep,\alpha} \rho_\alpha\bigr)(x-x') dx'.$$ 
We set 
$N_{\ep,\alpha}:=N_\ep * \rho_\alpha.$
In this way, we obtain a smooth term in $L^2(\R_+)$. Furthermore, $N_{\ep,\alpha}$   converges  to $N_\ep$ in $L^2(\R_+)$ when $\alpha$ tends to zero. 
We now have to find $B_{\ep,\alpha}$ solution of
\beq
\label{eq:inverse:filter}
4B_{\ep,\alpha}(2x)N_{\ep,\alpha}(2x)  - B_{\ep,\alpha}(x) N_{\ep,\alpha}(x)= \f{\p}{\p x} N_{\ep,\alpha} + \lb_{\ep,\alpha}  N_{\ep,\alpha}(x),\quad x\geq  0 \mbox{ .}
\eeq
Once again, the numerical aspects of the solution of such problem are addressed in Section~3 of \cite{DPZ}. 

\paragraph{Convergence and the Hybrid Method}
In \cite{PZ,DPZ} estimates are obtained that guarantee that for the choice
of 
$\alpha=\mathcal{O}(\sqrt{\epsilon})$ the computed value of $B_{\ep,\alpha}$ converges 
to $B$ in the appropriate norms, 
with a speed of convergence in the order of $\sqrt{\epsilon}$.  

In a number of numerical experiments with simulated data we found it convenient to
use a hybrid method where we perform filtering and apply the quasi-reversibility
technique. In the case of $g(x)=\kappa x$, at present we do not have the theoretical and numerical results of \cite{PZ,DPZ}. However, a simple adaptation of our codes led to implementations that yield
comparable results. We shall report on both cases.

\section{Reconstruction Results} \label{main}

In this section we report the results of the reconstructions in {\em two} 
sets of data from \cite{K69}. We actually performed the numerical inversion in
four sets, but we chose two of them where the experimental set up seemed
to fit better our model. 
As expected with real data, a number of difficulties have to be overcome. To cite a few:
a) The data is obtained on a rather scarce set of discrete points.
b) It clearly displays a substantial amount of noise. Yet, such noise is very
hard to quantify.
c)
We are not sure whether the conditions for the validity of the model were 
indeed satisfied, and whether enough time ellapsed so as to guarantee a good approximation
to stationary limiting distribution. 

\begin{figure}
\centering
\includegraphics [width=0.3\textwidth, height=4cm]{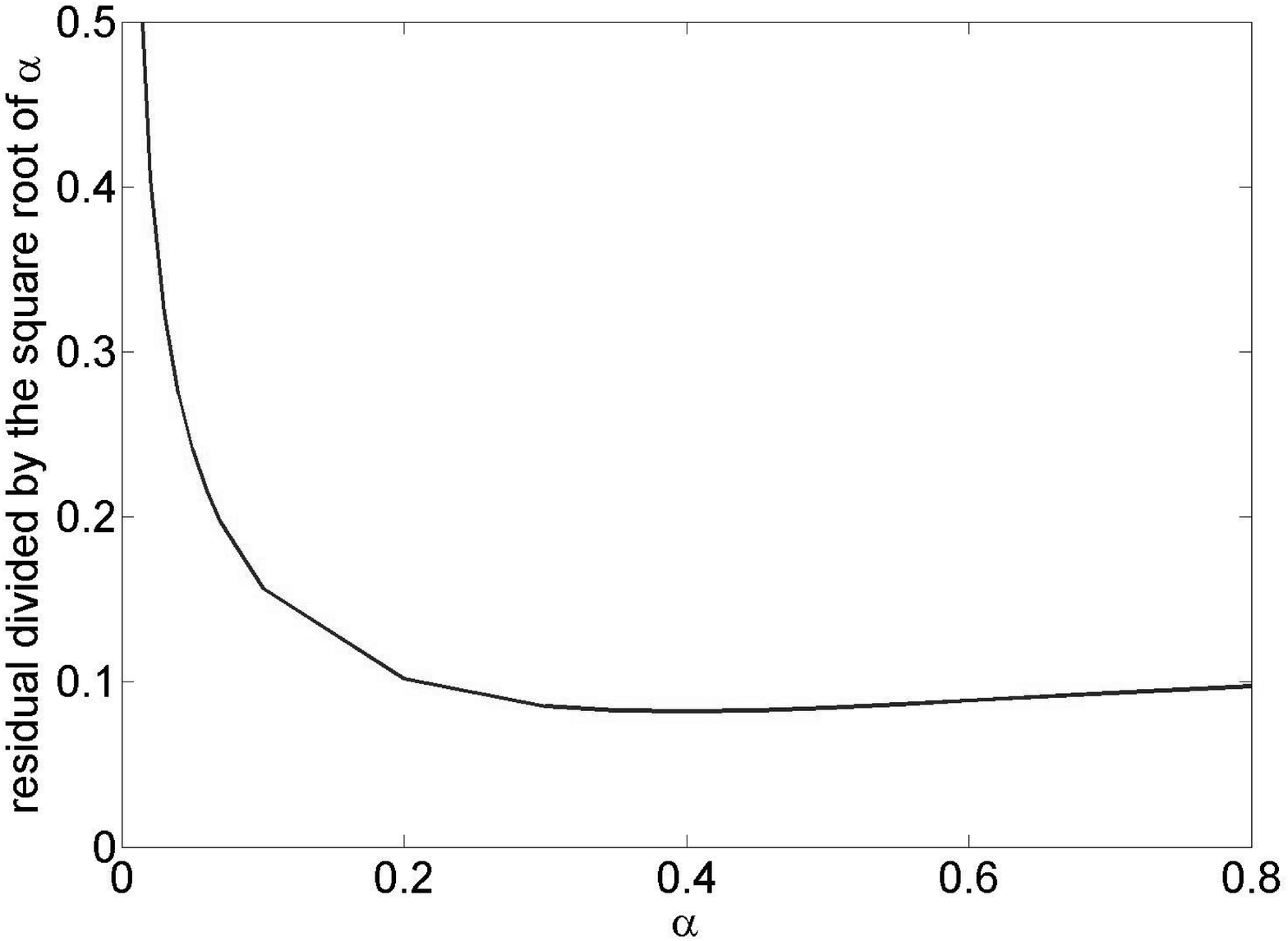}
\includegraphics [width=0.3\textwidth, height=4cm]{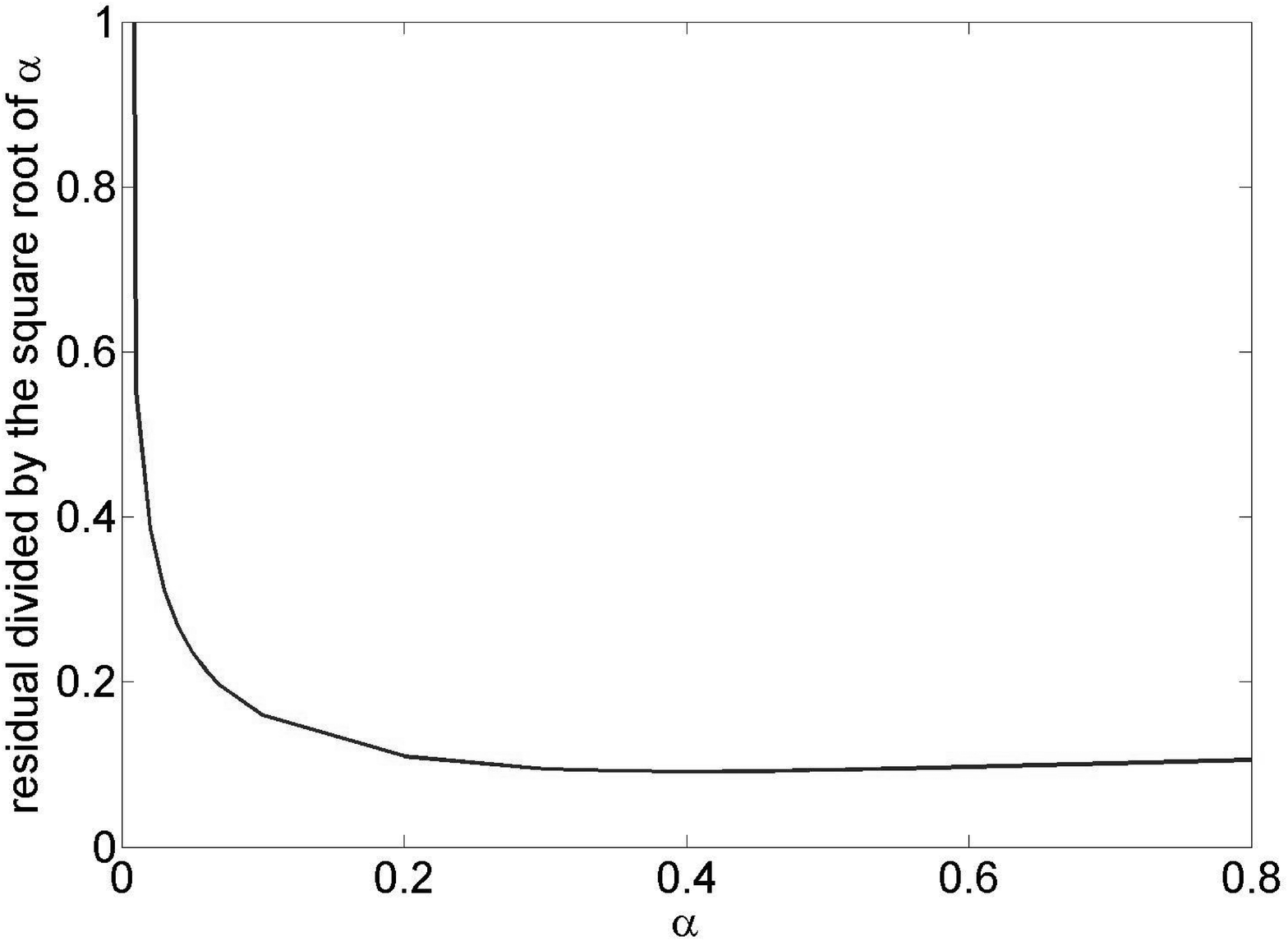}
\includegraphics [width=0.3\textwidth, height=4cm]{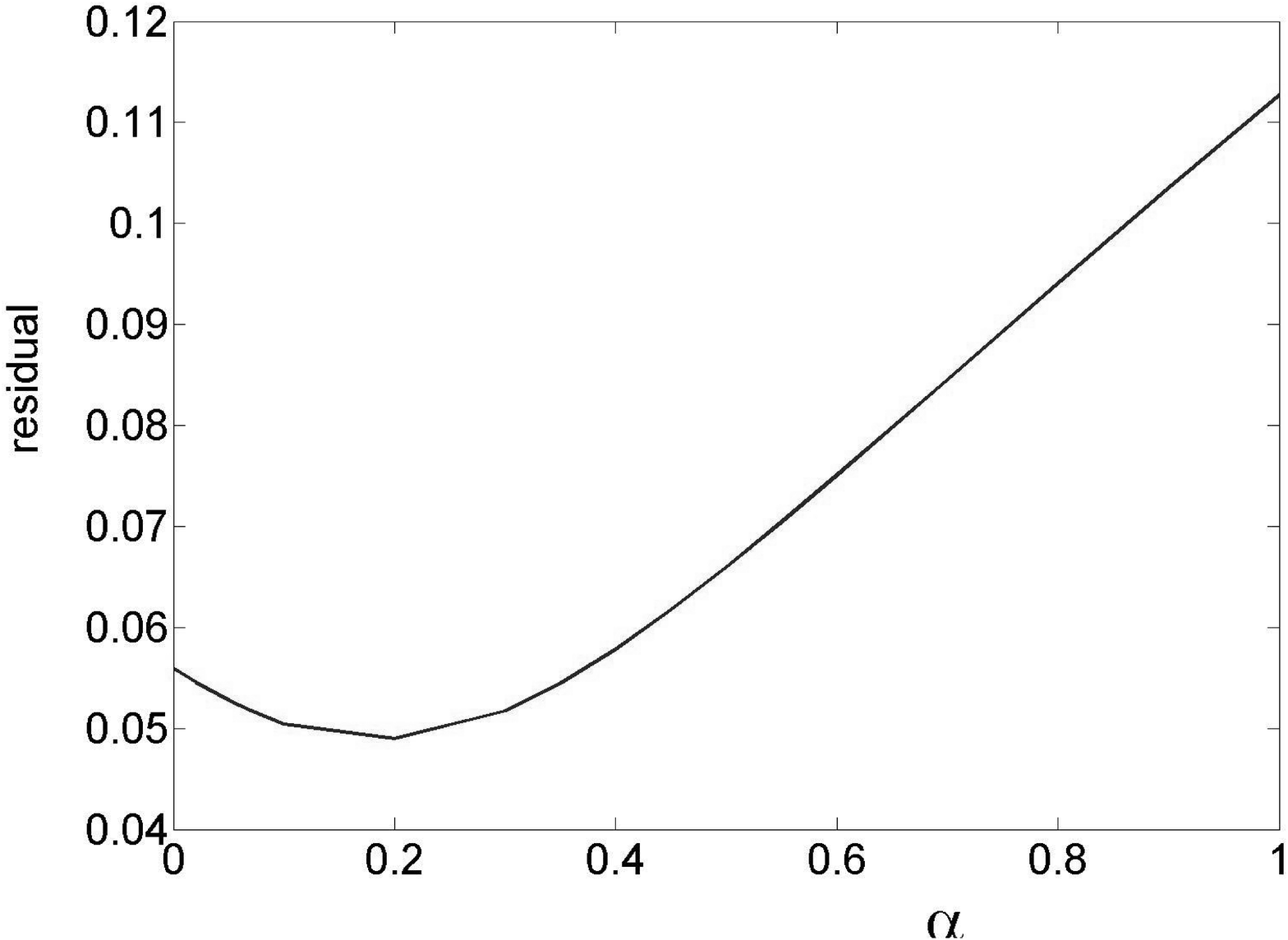}
\caption{\label{figLcurve} Example of criteria for the choice of $\alpha.$ Left: for a doubling time of 20 min, curve of $\f{residual}{\sqrt{\alpha}}$ with respect to $\alpha,$ linear case. The minimum is attained for $\alpha=0.4$ but the curve is very flat for $\alpha\geq 0.2.$ Center: same curve for doubling time of 54 min, linear case, the minimum is still for $\alpha=0.4.$ Right: as in Center, doubling time of 54 min and linear case, the curve of the residual exhibits a minimum for $\alpha=0.2.$}
 
\end{figure}

We performed the reconstructions both for the linear and the exponential models.
For each of the data sets, we followed the following steps (both for the linear model
as well as for the exponential model):
\begin{enumerate}
\item Transcribed the reported measurements from \cite{K69} to an array. 
\item Completed the boundary data for the volume $x$ close to zero and large enough
by setting the boundary conditions to zero (as required by the theory).
\item Interpolated to a uniform grid on the region under consideration using splines.
\item Deduced from \eqref{eq:lbn} and from the knowledge of both $N_\ep$ and $\lb_0$ (given in the article \cite{K69} by the doubling times $T_0 = \frac{Log(2)}{\lb_0}$) the constant value $g_0=\lb_0 \frac{\int x N dx}{\int N dx}$ in the definition of $g(x)=g_0$ (linear growth) or $g(x)=\lb_0 x$ (exponential growth).
\item Performed a search for a good regularization parameter $\alpha $ that would 
give a sufficiently small residual in the quasi-reversibility method of \cite{PZ} and in the hybrid methods. The result was also validated by other heuristical criteria
such as the use of 
the minimum of the ratio $\frac{residual}{\sqrt{\alpha}}$
\cite{EnHaNe}, see Figure \ref{figLcurve}.
\item Studied the behavior of the solutions to the inverse problem by varying the 
regularization parameter. 
\item Studied the consistency of the computed limiting distribution for the reconstructed
$B$ and the input data 
by computing the direct problem for some variants of the reconstructed $B.$
\end{enumerate}

We performed a number of extensive tests and examples, some of which that were not reported here could be found in 
\cite{Maia09}.
In all that follows we take the filter parameter to be 
$\alpha = 0.0001$. 
We remark that this choice is very arbitrary, but it is based on the standard deviation of the volumes once we assume the instrumental resolution of the Coulton counter has a standard deviation of the order $\sigma = 0.03 
\mu \mathrm{m}$ in the diameter \cite{K69}. Thus, the standard deviation of the volume,
assuming a normal distribution, is of the
order of $\sigma_V =(\pi/6) \sqrt{15} \sigma^3\approx 10^{-4} (\mu\mathrm{m})^3$.

%
%

\begin{figure}[h!]
\centering
\includegraphics[width=5cm, height=4cm]{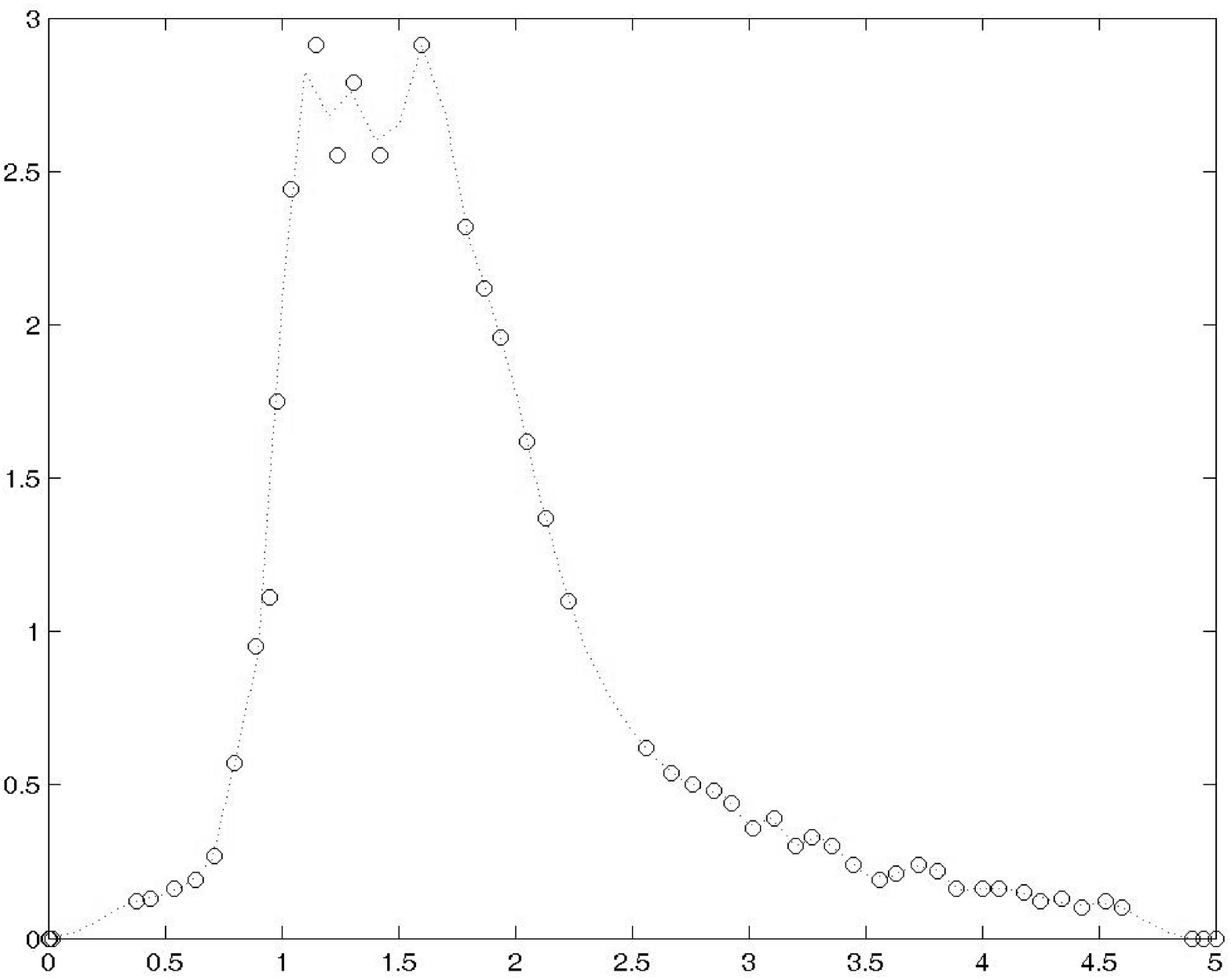}
\includegraphics[width=0.3\textwidth
, height=4cm]{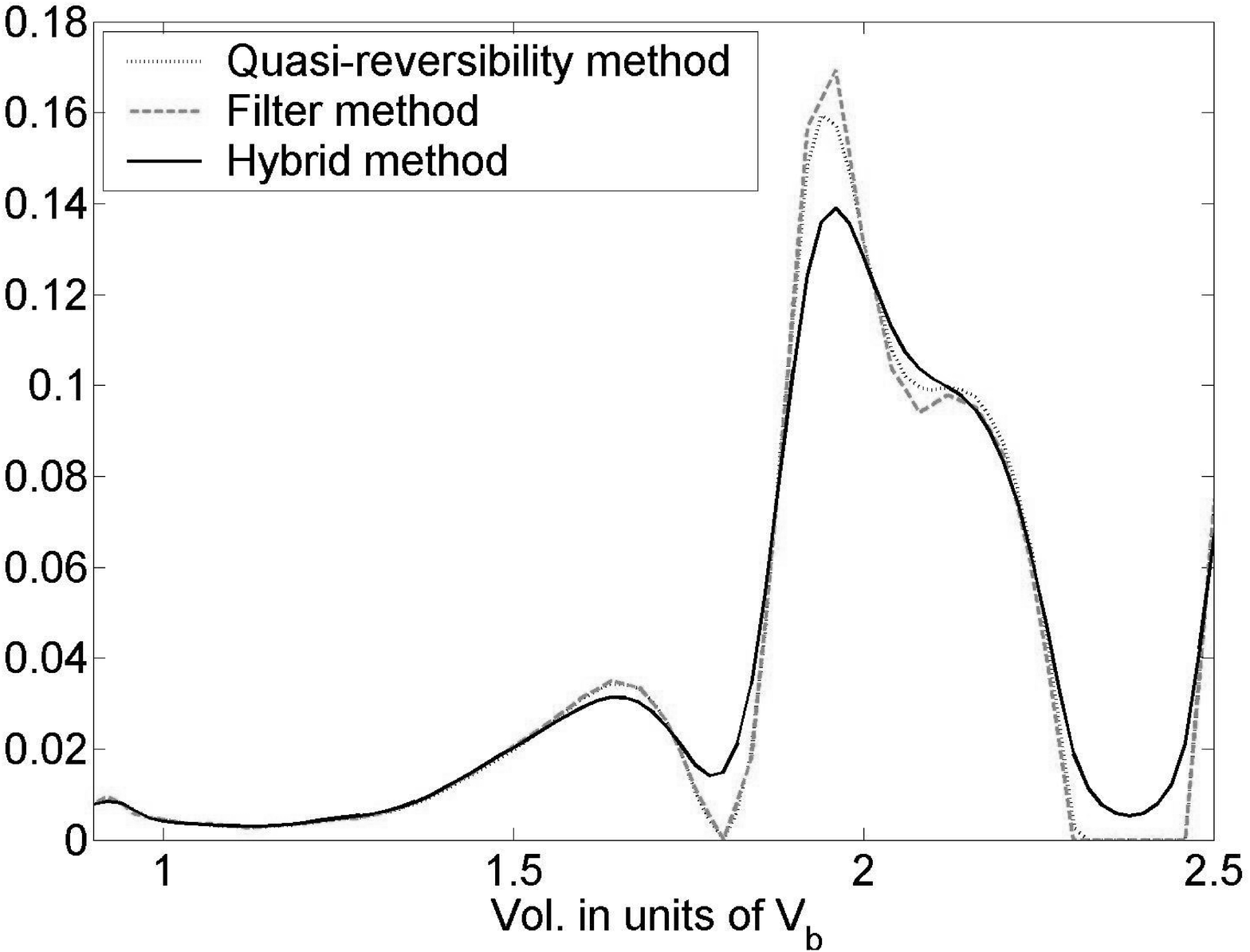}
\includegraphics [width=0.3\textwidth, 
height=4cm]{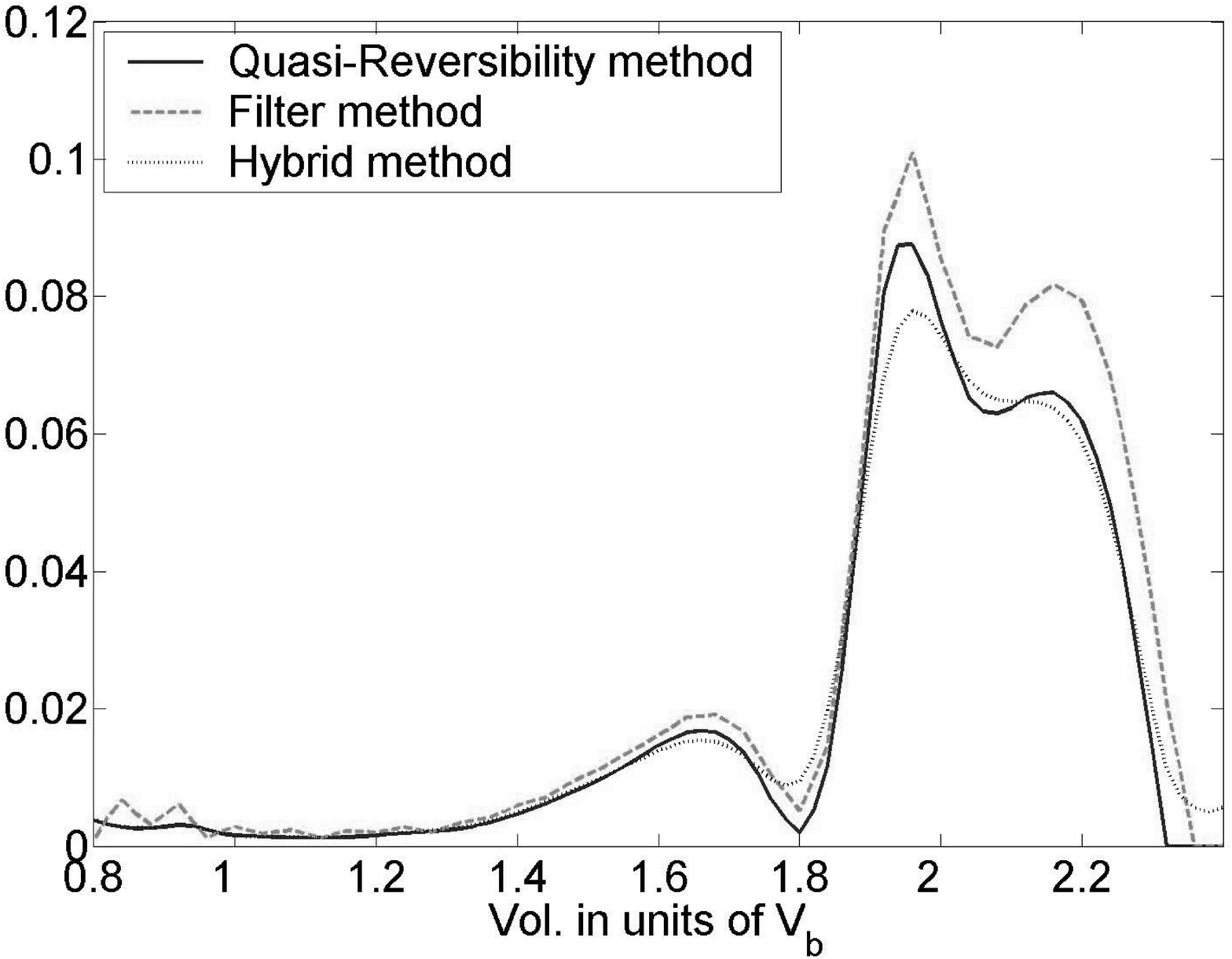}
\caption{\label{fig1T20}  Input data from~\cite{K69} for doubling time of 20 
min (left), reconstructed $B$ for the linear model (center), reconstructed $B$ for the exponential model (right).  $V_b = 1.36 (\mu \mathrm{m})^3.$ The value of 
$\alpha = 0.1$ (larger $\alpha$ give lower peaks).}
\end{figure}

\begin{figure}[h!]
\centering
\includegraphics [width=0.3\textwidth]{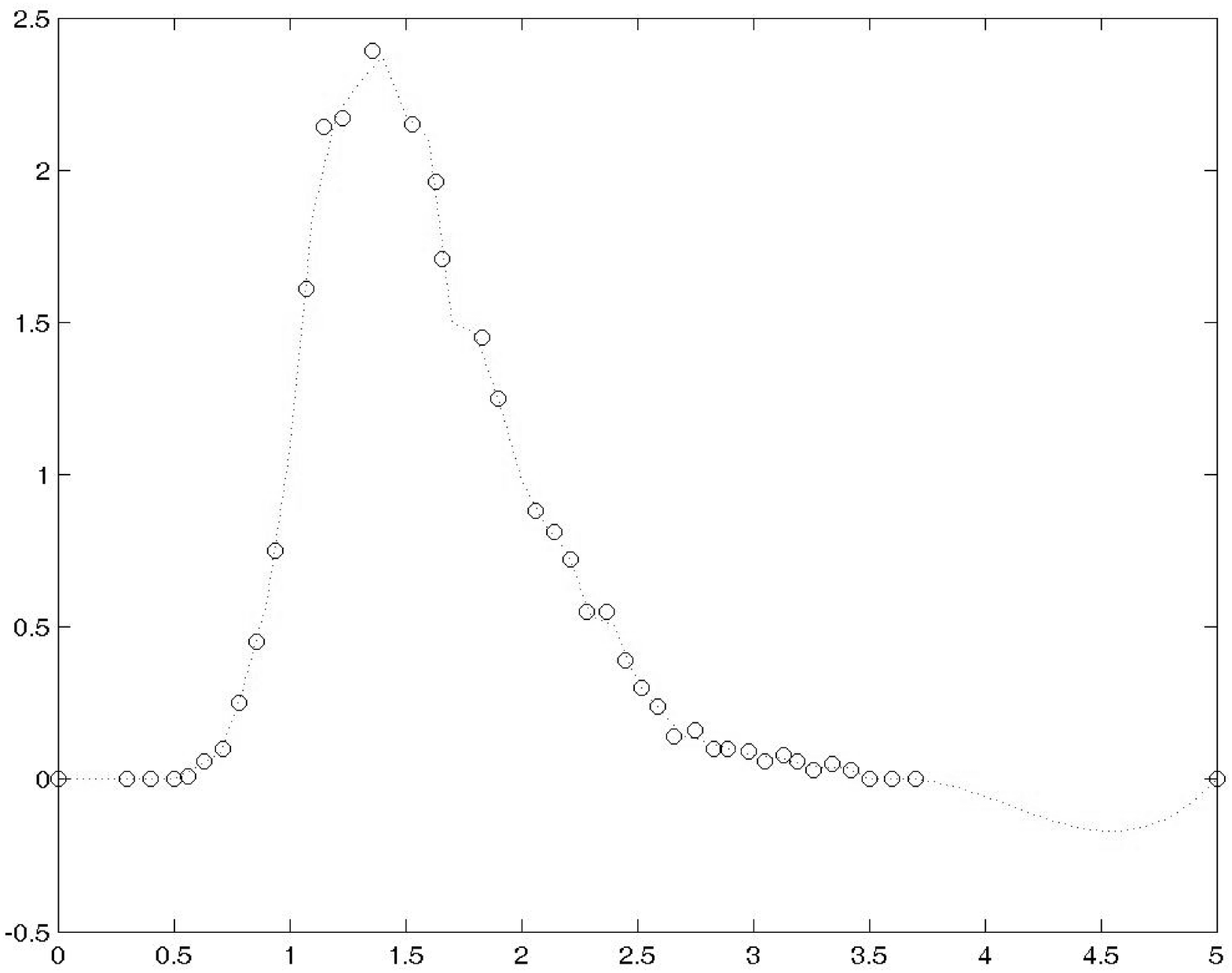}
\includegraphics [width=0.3\textwidth, 
height=4cm]{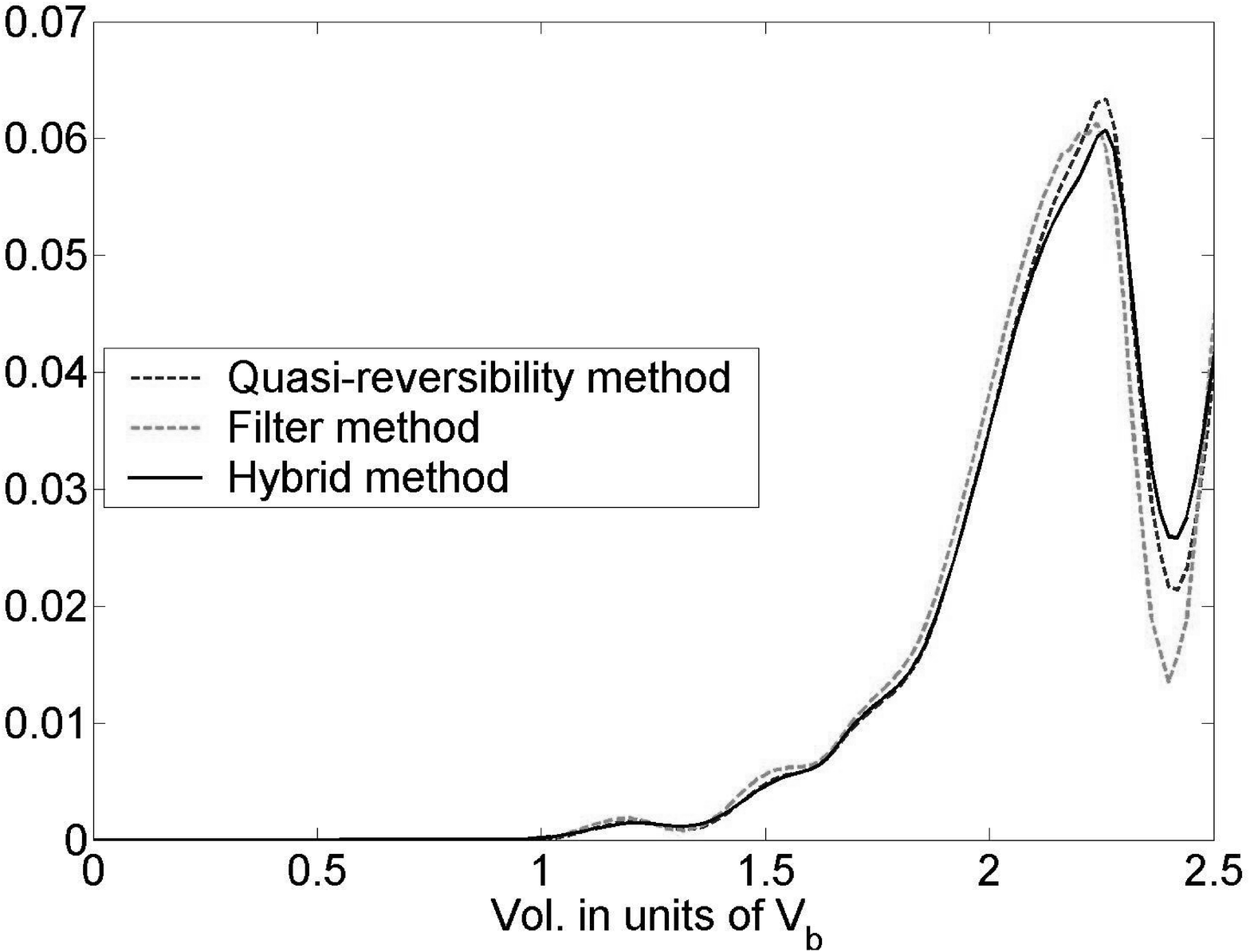}
\includegraphics [width=0.3\textwidth, 
height=4cm]{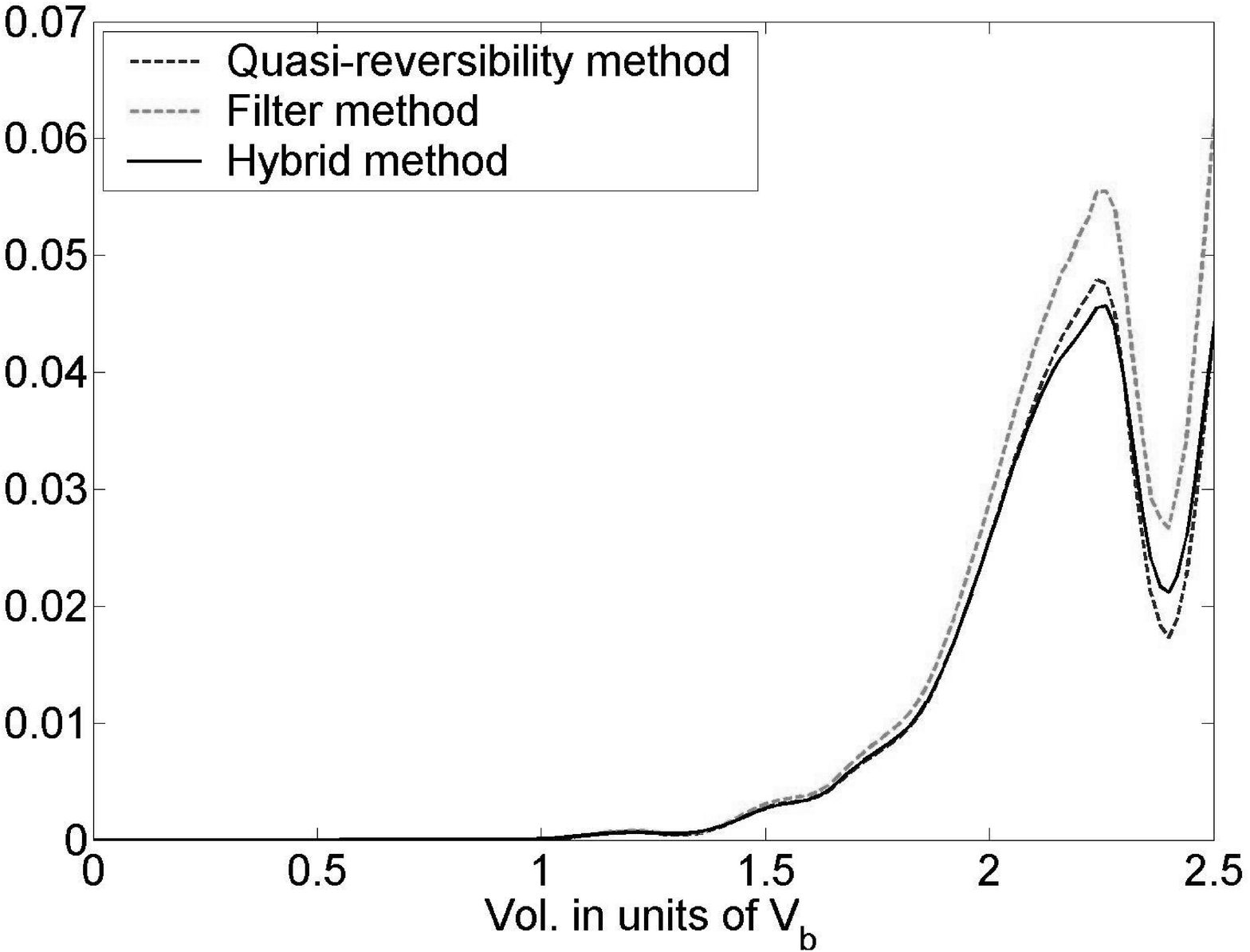}
\caption{\label{fig1T54} Input data from~\cite{K69} for doubling time of 54 min (left), reconstructed $B$ for the linear model, reconstructed $B$ for the exponential model (right). $V_b = 0.61 (\mu \mathrm{m})^3.$ The value of 
$\alpha = 0.2$ for the linear case, $\alpha=0.1$ for the exponential case.  }
\end{figure}


\section{Discussion and Perspectives} \label{disc} 

%
%

Most of the reconstructed birth rates we have obtained seem to
display a global maximum close to twice the average cell volume. 
Furthermore, in all cases we have at least a local maximum in that
neighborhood. The precise location, however, varies substantially with
the doubling time of the E. coli experiments, which in turn is a 
function of the experimental conditions. The peak locations also vary 
a little with the inversion methods and regularization parameters. 

The presence of doublets and triplets in the experimental data may give rise to
artifacts for values of $x$ higher than $3$. Indeed, it is not clear, due to the
highly nonlinear characteristic of the reconstruction whether an incorrect 
increase in the measurements near a certain $x_0$ would lead to an increase in the
reconstructed value of $B(x_0)$. 
This is a point that might be investigated further both theoretically and 
numerically.

We cannot infer any conclusion on the discussion of whether individual
cells undergo exponential or linear volume growth. Still, a superficial 
observation of the results seems to indicate that exponential growth gives smaller
residuals. This calls for a more extensive testing on better quality experimental
data. 

The optimal choice of the regularization parameter $\alpha$ should be investigated further. In this work we took a rather naive choice
due to lack of a good estimates on the experimental noise. Indeed, we tried to locate 
the optimal choice by looking at the curve of the residuals as a function of $\alpha$.
A natural extension of this work would be a more detailed investigation of {\em a priori}
and {\em a posteriori} choices for $\alpha$, 
such as the quasi-optimality criterion \cite{BaKi} or the L-curve method \cite{BaLe}.

\bibliographystyle{unsrt}
\bibliography{pronewzupedo}

%


\subsection*{Acknowledgments}
The authors were supported by the CNPq-INRIA agreement INVEBIO. 
JPZ was supported by CNPq under grants 302161/2003-1 and
474085/2003-1. JPZ  is thankful to the RICAM special semester
and to the International Cooperation Agreement Brazil-France.
A substantial part of this work was developed during a 3-month
international internship of PM at INRIA Rocquencourt during the Spring
of 2008 and supported by INRIA. PM and MD thank S. Boatto (UFRJ) for
facilitating this visit and for helpful discussions.

\end{document}